\documentclass[12pt,a4paper]{article}
\usepackage{amsmath,amssymb,amsthm}
\usepackage[
  pdftitle={On the periodic properties of self-decimated generators
of pseudorandom numbers},
  pdfauthor={Sergey Agievich, Oleg Solovey},
  bookmarks=true,pdfpagemode=None,pdfstartview=FitH]{hyperref}

\usepackage{eucal}

\newcommand{\expect}[1]{\mathop{\mbox{\normalfont\bfseries\sffamily E}}#1}
\newcommand{\disp}[1]{\mathop{\mbox{\normalfont\bfseries\sffamily D}}#1}
\newcommand{\PrSymbol}{\mbox{\sffamily\upshape\bfseries P}}
\newcommand{\prob}[2][]{\PrSymbol_{#1}\left\{ #2 \right\}}
\renewcommand{\vec}[1]{\mathbf{#1}}
\newcommand{\vecpi}{\boldsymbol{\pi}}

\makeatletter
\renewcommand{\@biblabel}[1]{#1.}
\makeatother

\newtheorem*{theorem*}{Theorem}

\begin{document}
\begin{sloppypar}

\title{On the periodic properties of self-decimated generators
of pseudorandom numbers}

\author{Sergey Agievich, Oleg Solovey\\
\small National Research Center
for Applied Problems of Mathematics and Informatics\\[-0.8ex]
\small Belarusian State University\\[-0.8ex]
\small Fr. Skorina av. 4, 220050 Minsk, Belarus\\[-0.8ex]
\small \texttt{agievich@bsu.by}, \texttt{solovey@bsu.by}
}

\date{}
\maketitle

\begin{abstract}
We consider a~self-decimated generator of pseudorandom numbers and
examine the preperiod $\lambda$ and the period $\mu$ of its state
sequence. We obtain the expectations and variances of $\lambda$
and $\mu$ for the case when decimation steps are chosen randomly
and independently from the set $\{1,2\}$.
\end{abstract}

\section{Results}\label{CLOCK.Results}

Let ${\mathbb N}$ be the set of all positive integers,
${\mathbb N}_0={\mathbb N}\cup\{0\}$,
$A\subset{\mathbb N}_0$ be some finite alphabet.
We denote by $A^*$ the set of all finite words
in the alphabet~$A$ (including the empty word~$\varepsilon$)
and by~$A^\omega$ the set of all one-way infinite words
(see~\cite{GouJac1983,ChoKar1997} for further details).
Given a word~$a\in A^*$,
let~$l(a)$ be the length of~$a$ and
$w(a)$ be the sum of letters of~$a$.
Let~$c^n$ denote~$n$ successive instances of the letter~$c$
and~$ab$ denote the concatenation of the words~$a$ and~$b$.

Given~$s=s_0 s_1\ldots\in A^\omega$ and $T\in\mathbb{N}$,
determine the numbers $\lambda=\lambda(s,T)\in{\mathbb N}_0$ and
$\mu=\mu(s,T)\in{\mathbb N}$ such that
\begin{equation}\label{Eq.CLOCK.1}
T\mid(s_{\lambda}+s_{\lambda+1}+\ldots+s_{\lambda+\mu-1})
\end{equation}
and $T$ does not divide any of the sums
$s_t+s_{t+1}+\ldots+s_{\tau}$, $0\leq t<\tau<\lambda+\mu-1$.
For example, if $s=2212221\ldots=2^2 1 2^3 1\ldots$, then
$$
\lambda(s,1)=0,\ \mu(s,1)=1,\quad
\lambda(s,2)=0,\ \mu(s,2)=1,\ldots\quad
\lambda(s,8)=2,\ \mu(s,8)=5,\ldots.
$$

The characteristics~$\lambda$ and $\mu$ describe
periodic properties of self-decimated generators
of pseudorandom numbers~(see~\cite{ChaGol1998,Rue1988}).
Consider a generator~$G$ with the set of internal states
${\cal S}$, $|{\cal S}|=T$, and the state-transition function
$\varphi\colon{\cal S}\to{\cal S}$.
Let $\varphi$ be a full cycle substitution, i.e.,
for any~$S\in{\cal S}$ the images
$\varphi(S)$, $\varphi^2(S)=\varphi(\varphi(S)),\ldots$, $\varphi^T(S)$
run over all ${\cal S}$.
Under the usual way of $G$~functioning,
one chooses an initial state~$S_0\in{\cal S}$,
calculates the sequence
\begin{equation}\label{Eq.CLOCK.G1}
S_{t+1}=\varphi(S_t),\quad t=0,1,\ldots,
\end{equation}
and uses the current internal state~$S_t$ to determine
the current output pseudorandom number.
The {\it self-decimation} (of internal states of $G$) means
using an additional function~$d\colon{\cal S}\to A$
and replacing~\eqref{Eq.CLOCK.G1} by the rule
\begin{equation}\label{Eq.CLOCK.G2}
S_{t+1}=\varphi^{d(S_t)}(S_t),\quad t=0,1,\ldots.
\end{equation}
Now, if the word~$s$ consists of the successive letters~$s_t=d(S_t)$,
then $\lambda(s,T)$ and $\mu(s,T)$ are respectively
the preperiod and period of the state sequence~\eqref{Eq.CLOCK.G2}.

If the letters of~$s$ are chosen randomly,
then~$\lambda$ and~$\mu$ become random variables.
The expectations~$\expect\lambda$, $\expect\mu$
and variances~$\disp\lambda$, $\disp\mu$
of these variables are of interest
in connection with the estimation of the
preperiod and period of self-decimated generators.
One of the well-known results in the theory of random mappings can
be stated as follows: If the letters of~$s$ are chosen randomly,
independently and uniformly from the
alphabet~$A=\{0,1,\ldots,T-1\}$, then $\expect\lambda(s,T)$
and~$\expect\mu(s,T)$ have the form $\sqrt{\pi T/8}+O(1)$
as~$T\to\infty$ (see, for example,~\cite{FlaOdl1990}).
Unfortunately, obtaining similar asymototic (or exact)
expressions for arbitrary alphabet~$A$
seems to be a more complex task.

We consider the frequent choice $A=\{1,2\}$ and obtain
the following result.

\begin{theorem*}
If the letters of the word $s=s_0s_1\ldots\in \{1,2\}^\omega$
are chosen randomly, independently with the probabilities
$$
\prob{s_t=1}=1-p,\quad
\prob{s_t=2}=p,\quad
0<p<1,\quad
t=0,1,\ldots,
$$
then as~$T\to\infty$ it holds that
\begin{align}
\label{Eq.CLOCK.Expect}
\expect\lambda(s,T)&=\frac{p}{(1-p)(1+p)^2}+O(p^TT),&
%
%
%
\expect\mu(s,T)&=\frac{T}{1+p}+O(p^TT),\\
%
%
\label{Eq.CLOCK.Disp}
\disp\lambda(s,T)&=\frac{p+p^3+p^4}{(1-p)^2(1+p)^4}+O(p^T T^2),&
\disp\mu(s,T)&=\frac{Tp(1-p)}{(1+p)^3}+O(p^T T^2).
\end{align}
\end{theorem*}

Note that~\eqref{Eq.CLOCK.Expect}, \eqref{Eq.CLOCK.Disp}
also hold for the random word $s=s_0s_1\ldots\in\{q,2q\}^\omega$,
where~$q\in{\mathbb N}$ is coprime to~$T$ and
$\prob{s_t=q}=1-p$, $\prob{s_t=2q}=p$.

Note also that R.~Rueppel in the paper~\cite{Rue1988}
examined~the case when the word~$s$ is determined
from a linear recurrence sequence~$\sigma_0,\sigma_1,\ldots$
over the field of order~$2$ by the rule
$$
s_t=
\begin{cases}
1, & \sigma_t=0,\\
2, & \sigma_t=1,
\end{cases}
\quad t=0, 1,\ldots.
$$
Rueppel showed that if the characteristic polynomial of the
sequence $\sigma_0,\sigma_1,\ldots$ is primitive of degree $k$,
then
$$
\mu(s,T)=\left\lfloor\frac{2T}{3}\right\rfloor,\quad
T=2^k-1,
$$
where  $\lfloor z\rfloor$ is the largest integer~$\leq z$.
As we can see from Theorem,
this estimation agrees with the expectation~$\expect\mu(s,T)$
for the random word $s$ with~$p=1/2$.

\section{Proof}\label{CLOCK.Proof}

Given~$s$ and~$T$, determine
the words~$s_0\ldots s_{\lambda-1}$
and~$s_{\lambda}\ldots s_{\lambda+\mu-1}$,
which we call the {\it prefix} and {\it cyclic part} of $s$,
respectively.
Conversely, any nonempty word~$a\in A^*$ is a cyclic part
of some word $s$.
Given~$a$, we can determine the possible values of $T=T(a)$ and the
set~$B(a)$ of the possible prefixes of~$s$,
using the following restrictions:
\begin{itemize}
\item[P1)]
$w(a)\equiv 0\pmod{T}$;

\item[P2)]
the residues
\begin{equation}\label{Eq.CLOCK.GF.1}
0,\quad a_0,\quad a_0+a_1,\ldots,\quad
a_0+a_1+\ldots+a_{l(a)-2}\pmod{T}
\end{equation}
are pairwise distinct;

\item[P3)]
if $b\in B(a)$ and $b\neq\varepsilon$, then any of the residues
$$
-b_{l(b)-1},\quad -b_{l(b)-1}-b_{l(b)-2},\ldots,\quad
-b_{l(b)-1}-b_{l(b)-2}-\ldots-b_0\pmod{T}
$$
differs from the residues~\eqref{Eq.CLOCK.GF.1}.
\end{itemize}

Indeed, the condition P1 follows from the definition of a cyclic part.
If the condition P2 fails, then
$$
a_t+\ldots+a_\tau\equiv 0\pmod{T}
$$
for some $t$, $\tau$, $0\leq t<\tau<l(a)-1$
and $a$ cannot be a~cyclic part.
Finally, if the condition P3 fails,
then one of the congruences
$$
b_t+\ldots+b_{l(b)-1}\equiv 0\pmod{T},\quad
b_t+\ldots+b_{l(b)-1}+a_0+\ldots+a_\tau\equiv 0\pmod{T}
$$
holds for some~$t\leq l(b)-1$, $\tau<l(a)-1$
and~$a$ cannot be a~cyclic part again.

Suppose that $A=\{1,2\}$ and let $\vec{x}=(x_1, x_2)$, $\vec{y}=(y_1, y_2)$.
Given $\Omega\subseteq A^*$, define the generating function
\begin{equation}\label{Eq.CLOCK.G}
G_\Omega(\vec{x},\vec{y},z)=
\sum_{\scriptstyle n_1,n_2,m_1,m_2\geq 0\atop\scriptstyle t\geq 1}
g(n_1,n_2,m_1,m_2,t) x_1^{n_1} x_2^{n_2} y_1^{m_1} y_2^{m_2} z^t,
\end{equation}
where $g(n_1,n_2,m_1,m_2,t)$ is the number of the words $s$ with
the cyclic part $a\in\Omega$ and prefix $b\in B(a)$
such that $a$ contains $n_1$ letters~$1$ and $n_2$ letters~$2$,
$b$ contains $m_1$ letters~$1$ and $m_2$ letters~$2$,
and~$T(a)=t$.

Let us divide $A^*$ into the subsets $\Omega_1$, $\Omega_2$, $\Omega_3$,
$\Omega_4$ (which will be defined below)
and determine the generating function of the form~\eqref{Eq.CLOCK.G}
for each subset.

\begin{enumerate}
\item
$\Omega_1=\{2\}^*1$.
If~$a\in\Omega_1$, then $T\mid w(a)$ according to
P1 and $T\geq l(a)>w(a)/2$ according to P2.
Therefore, $T(a)=w(a)$.
Now, if $a=2^m 1$, then
$B(a)=\{2^k\colon k=0,\ldots,m\}$ according to P3.
The target generating function has the form
\begin{align*}
G_{\Omega_1}(\vec{x},\vec{y},z)&=
\sum_{m\geq 0} x_1 x_2^{m} z^{2m+1}\sum_{k=0}^{m} y_2^k\\
&=\frac{x_1z}{1-y_2}\sum_{m\geq 0}(x_2z^2)^m-
\frac{x_1y_2z}{1-y_2}\sum_{m\geq 0}(x_2y_2z^2)^m\\
&=\frac{x_1 z}{(1-y_2)(1-x_2 z^2)}-\frac{x_1 y_2 z}{(1-y_2)(1-x_2 y_2 z^2)}\\
&=\frac{x_1 z}{(1-x_2 z^2)(1-x_2 y_2 z^2)}.
\end{align*}

\item
$\Omega_2=\{1,2\}^*1\{2\}^*1$.
Again, $T(a)=w(a)$ according to P1, P2.
Let $a=\alpha 12^m1$, $\alpha\in A^*$.
Using P3, we obtain~$B(a)=\{2^k\colon k=0,\ldots,m\}$.
Therefore,
\begin{align*}
G_{\Omega_2}(\vec{x},\vec{y},z)
&=\sum_{l\geq 0}(x_1 z + x_2 z^2)^l
\sum_{m\geq 0}x_1^2 x_2^{m}z^{2m+2}\sum_{k=0}^{m}y_2^k\\
&=\frac{x_1^2 z^2}{(1-x_1z-x_2z^2)(1-x_2z^2)(1-x_2y_2z^2)}.
\end{align*}

\item $\Omega_3=\{1,2\}^*12\{2\}^*$. In this case $T(a)=w(a)$ and
$B(a)=\{\varepsilon\}\cup\{2^k1\colon k=0,\ldots,m\}$
for~$a=\alpha 122^m$.
Therefore,
\begin{align*}
G_{\Omega_3}(\vec{x},\vec{y},z)&=
\sum_{l\geq 0}(x_1z+x_2z^2)^l
\sum_{m\geq 0}x_1 x_2^{m+1}z^{2m+3} \left(1+\sum_{k=0}^{m}y_1 y_2^k \right)\\
&=\frac{x_1x_2z^3(1+y_1-x_2y_2z^2)}{(1-x_1z-x_2z^2)(1-x_2z^2)(1-x_2y_2z^2)}.
\end{align*}

\item $\Omega_4=\{2\}^*$. If $a=2^m$, then $T(a)=2m$ and
$B(a)=\{\varepsilon\}\cup\{2^k1\colon k=0,\ldots,m-1\}$.
In addition, $T(a)=m$ and $B(a)=\{\varepsilon\}$ for odd~$m$.
The generating function has the form
\begin{align*}
G_{\Omega_4}(\vec{x},\vec{y}, z) &=\sum_{m\geq1}
x_1^{m}z^{2m} \left(1+\sum_{k=0}^{m-1} y_1 y_2^k \right)+
\sum_{m\geq 0}x_2^{2m+1}z^{2m+1}\\
&=\frac{1+x_2y_1z^2-x_2y_2z^2}{(1-x_2z^2)(1-x_2y_2z^2)}+
\frac{x_2z}{1-x_2^2z^2}.
\end{align*}
\end{enumerate}

Finally, we get
\begin{gather*}
G_{A^*}(\vec{x},\vec{y},z)
=G_{\Omega_1}(\vec{x},\vec{y},z)+ G_{\Omega_2}(\vec{x},\vec{y},z)+
 G_{\Omega_3}(\vec{x},\vec{y},z)+ G_{\Omega_4}(\vec{x},\vec{y},z)\\
=\frac{1+x_2y_1z^2+x_1x_2y_2z^3-x_2y_2z^2}{(1-x_2y_2z^2)(1-x_1z-x_2z^2)}
+ \frac{x_2z}{1-x_2^2z^2}.
\end{gather*}

To find the expectation $\expect\lambda(s,T)$,
introduce the operator
$E_{\vec{y}}=
y_1\frac{\partial}{\partial y_1}+y_2\frac{\partial}{\partial y_2}$
and denote $\vecpi=(1-p,p)$.
We have
\begin{gather*}
\expect\lambda(s,T)=\sum_{n_1, n_2, m_1, m_2 \geq 0}
\frac{g(n_1,n_2,m_1,m_2,T)(m_1+m_2)}{(1-p)^{n_1+m_1}p^{n_2+m_2}}=
[z^T]\left.E_{\vec{y}}G_{A^*}(\vecpi,\vec{y},z)
\right|_{\vec{y}=\vecpi},
\end{gather*}
where $[z^T]f(z)$ is the coefficient of $z^T$ in
$f(z)=f_0+f_1 z +f_2 z^2+\ldots$.
Now the first part of~\eqref{Eq.CLOCK.Expect}
follows from
\begin{align*}
\left.E_{\vec{y}}G_{A^*}(\vecpi,\vec{y},z)
\right|_{\vec{y}=\vecpi}&=
\frac{(1-p)pz^2}{(1-z)(1-p^2z^2)^2}\\
&=\frac{p}{(1-p)(1+p)^2(1-z)}-\frac{1}{4(1-pz)^2}-\frac{p}{4(1-p)(1-pz)}\\
&+\frac{1-p}{4(1+p)(1+pz)^2}-\frac{(1-p)p}{4(1+p)^2(1+pz)}\\
&=\frac{p}{(1-p)(1+p)^2}\sum_{k\geq 0}z^k-
\frac{1}{4}\sum_{k\geq 0}(k+1)p^kz^k-
\frac{p}{4(1-p)}\sum_{k\geq 0}p^kz^k\\
&+\frac{1-p}{4(1+p)}\sum_{k\geq 0}(k+1)(-p)^kz^k-
\frac{(1-p)p}{4(1+p)^2}\sum_{k\geq 0}(-p)^kz^k.
\end{align*}

Similarly, introducing the operator
$E_{\vec{x}}=
x_1\frac{\partial}{\partial x_1}+x_2\frac{\partial}{\partial x_2}$
and using the equalities
\begin{align*}
\expect\mu(s,T)&=[z^T]E_{\vec{x}}
\left.G_{A^*}(\vec{x},\vecpi,z)\right|_{\vec{x}=\vecpi},\\
\disp\lambda(s,T)&=[z^T]E^2_{\vec{y}}
\left.G_{A^*}(\vecpi,\vec{y},z)\right|_{\vec{y}=\vecpi}-
\left(\expect\lambda(s,T)\right)^2,\\
\disp\mu(s,T)&=[z^T]E^2_{\vec{x}}
\left.G_{A^*}(\vec{x},\vecpi,z)\right|_{\vec{x}=\vecpi}-
\left(\expect\mu(s,T)\right)^2,
\end{align*}
we prove the remaining parts of~\eqref{Eq.CLOCK.Expect},
\eqref{Eq.CLOCK.Disp}.

\end{sloppypar}
\end{document}